\documentclass[11pt]{article}
\usepackage{amsmath}
\usepackage{amssymb}
\usepackage{amsthm}
\numberwithin{equation}{section}
\begin{document}
\title{Riesz transforms associated with Schr\"odinger operators acting on weighted Hardy spaces}
\author{Hua Wang\,\footnote{E-mail address: wanghua@pku.edu.cn.}\\
\footnotesize{School of Mathematical Sciences, Peking University, Beijing 100871, China}}
\date{}
\maketitle
\begin{abstract}
Let $L=-\Delta+V$ be a Schr\"odinger operator acting on $L^2(\mathbb R^n)$, $n\ge1$, where $V\not\equiv 0$ is a nonnegative locally integrable function on $\mathbb R^n$. In this article, we will introduce weighted Hardy spaces $H^p_L(w)$ associated with $L$ by means of the area integral function and study their atomic decomposition theory. We also show that the Riesz transform $\nabla L^{-1/2}$ associated with $L$ is bounded from our new space $H^p_L(w)$ to the classical weighted Hardy space $H^p(w)$ when $\frac{n}{n+1}<p<1$ and$w\in A_1\cap RH_{(2/p)'}$.
\end{abstract}

\section{Introduction}
Let $n\ge1$ and $V$ be a nonnegative locally integrable function defined on $\mathbb R^n$, not identically zero. We define the form $\mathcal Q$ by
\begin{equation*}
\mathcal Q(u,v)=\int_{\mathbb R^n}\nabla u\cdot\nabla v\,dx+\int_{\mathbb R^n}Vuv\,dx
\end{equation*}
with domain $\mathcal D(\mathcal Q)=\mathcal V\times\mathcal V$ where
\begin{equation*}
\mathcal V=\{u\in L^2(\mathbb R^n):\frac{\partial u}{\partial x_k}\in L^2(\mathbb R^n)\mbox{ for }k=1,\ldots,n \mbox{ and }\sqrt {V}u\in L^2(\mathbb R^n)\}.
\end{equation*}
It is well known that this symmetric form is closed. Note also that it was shown by Simon [17] that this form coincides with the minimal closure of the form given by the same expression but defined on $C^\infty_0(\mathbb R^n)$(the space of $C^\infty$ functions with compact supports). In other words, $C^\infty_0(\mathbb R^n)$ is a core of the form $\mathcal Q$.

Let us denote by $L$ the self-adjoint operator associated with $\mathcal Q$. The domain of $L$ is given by
\begin{equation*}
\mathcal D(L)=\{u\in \mathcal D(\mathcal Q):\exists \,v\in L^2\mbox{ such that }\mathcal Q(u,\varphi)=\int_{\mathbb R^n}v\varphi\,dx,\forall\varphi\in\mathcal D(\mathcal Q)\}.
\end{equation*}
Formally, we write $L=-\Delta+V$ as a Schr\"odinger operator with potential $V$. Let $\{e^{-tL}\}_{t>0}$ be the semigroup of linear operators generated by $-L$ and $p_t(x,y)$ be their kernels. Since $V$ is nonnegative, the Feynman-Kac formula implies that
\begin{equation}
0\le p_t(x,y)\le\frac{1}{(4\pi t)^{n/2}}e^{-\frac{|x-y|^2}{4t}}
\end{equation}
for all $t>0$ and $x,y\in\mathbb R^n$.

The operator $\nabla L^{-1/2}$ is called the Riesz transform associated with $L$, which is defined by
\begin{equation}
\nabla L^{-1/2}(f)(x)=\frac{1}{\sqrt\pi}\int_0^\infty \nabla e^{-tL}(f)(x)\frac{dt}{\sqrt t}.
\end{equation}
This operator is bounded on $L^2(\mathbb R^n)$(see [11]). Moreover, it was proved in [1,3] that by using the molecular decomposition of functions in the Hardy space $H^1_L(\mathbb R^n)$, the operator $\nabla L^{-1/2}$ is bounded from $H^1_L(\mathbb R^n)$ into $L^1(\mathbb R^n)$, and hence, by interpolation, is bounded on $L^p(\mathbb R^n)$ for all $1<p\le2$. Now assume that $V\in RH_q$(Reverse H\"older class). In [15], Shen showed that $\nabla L^{-1/2}$ is a Calder\'on-Zygmund operator if $q\ge n$. When $\frac n2\le q<n$, $\nabla L^{-1/2}$ is bounded on $L^p(\mathbb R^n)$ for $1<p\le p_0$, where $1/{p_0}=1/q-1/n$, and the above range of $p$ is optimal. For more information about the Hardy spaces $H^p_L(\mathbb R^n)$ associated with Schr\"odinger operators for $0<p\le1$, we refer the readers to [4,5,6].

In [18], Song and Yan introduced the weighted Hardy spaces $H^1_L(w)$ associated to $L$ in terms of the area integral function and established their atomic decomposition theory. In the meantime, they also showed that the Riesz transform $\nabla L^{-1/2}$ is bounded on $L^p(w)$ for $1<p<2$, and bounded from $H^1_L(w)$ to the classical weighted Hardy space $H^1(w)$.

As a continuation of [18], the main purpose of this paper is to define the weighted Hardy spaces $H^p_L(w)$ associated to $L$ for $0<p<1$ and study their atomic characterizations. We also obtain that $\nabla L^{-1/2}$ is bounded from $H^p_L(w)$ to the classical weighted Hardy space $H^p(w)$ for $\frac{n}{n+1}<p<1$. Our main result is stated as follows.
\newtheorem{theorem}{Theorem}[section]
\begin{theorem}
Suppose that $L=-\Delta+V$. Let $\frac{n}{n+1}<p<1$ and $w\in A_1\cap RH_{(2/p)'}$. Then the operator $\nabla L^{-1/2}$ is bounded from $H^p_L(w)$ to the classical weighted Hardy space $H^p(w)$.
\end{theorem}
It is worth pointing out that when $L=-\Delta$ is the Laplace operator on $\mathbb R^n$, then the space $H^p_L(w)$ coincides with the classical weighted Hardy space $H^p(w)$. Therefore, in this particular case, we derive that the classical Riesz transform $\nabla(-\Delta)^{-1/2}$ is bounded on $H^p(w)$ for $\frac{n}{n+1}<p<1$, which was already obtained by Lee and Lin in [12].

\section{Notations and preliminaries}
First, let us recall some standard definitions and notations. The classical $A_p$ weight
theory was first introduced by Muckenhoupt in the study of weighted
$L^p$ boundedness of Hardy-Littlewood maximal functions in [13].
A weight $w$ is a locally integrable function on $\mathbb R^n$ which takes values in $(0,\infty)$ almost everywhere, $B=B(x_0,r)$ denotes the ball with the center $x_0$ and radius $r$.
We say that $w\in A_p$,\,$1<p<\infty$, if
$$\left(\frac1{|B|}\int_B w(x)\,dx\right)\left(\frac1{|B|}\int_B w(x)^{-\frac{1}{p-1}}\,dx\right)^{p-1}\le C \quad\mbox{for every ball}\; B\subseteq \mathbb
R^n,$$ where $C$ is a positive constant which is independent of $B$.

For the case $p=1$, $w\in A_1$, if
$$\frac1{|B|}\int_B w(x)\,dx\le C\cdot\underset{x\in B}{\mbox{ess\,inf}}\,w(x)\quad\mbox{for every ball}\;B\subseteq\mathbb R^n.$$

A weight function $w$ is said to belong to the reverse H\"{o}lder class $RH_r$ if there exist two constants $r>1$ and $C>0$ such that the following reverse H\"{o}lder inequality holds
$$\left(\frac{1}{|B|}\int_B w(x)^r\,dx\right)^{1/r}\le C\left(\frac{1}{|B|}\int_B w(x)\,dx\right)\quad\mbox{for every ball}\; B\subseteq \mathbb R^n.$$

If $w\in A_p$ with $1<p<\infty$, then we have $w\in A_r$ for all $r>p$, and $w\in A_q$ for some $1<q<p.$ It follows from H\"{o}lder's inequality that $w\in RH_r$ implies $w\in RH_s$ for all $1<s<r$. Moreover, if $w\in RH_r$, $r>1$, then we have $w\in RH_{r+\varepsilon}$ for some $\varepsilon>0$.

Given a ball $B$ and $\lambda>0$, $\lambda B$ denotes the ball with the same center as $B$ whose radius is $\lambda$ times that of $B$. For a given weight function $w$, we denote the Lebesgue measure of $B$ by $|B|$ and the weighted measure of $B$ by $w(B)$, where $w(B)=\int_Bw(x)\,dx.$

We give the following results that we will use in the sequel.

\newtheorem{lemma}[theorem]{Lemma}
\begin{lemma}[{[8]}]
Let $w\in A_p$, $p\ge1$. Then, for any ball $B$, there exists an absolute constant $C$ such that $$w(2B)\le C\,w(B).$$ In general, for any $\lambda>1$, we have $$w(\lambda B)\le C\cdot\lambda^{np}w(B),$$where $C$ does not depend on $B$ nor on $\lambda$.
\end{lemma}
\begin{lemma}[{[8,9]}]
Let $w\in A_p\cap RH_r$, $p\ge1$ and $r>1$. Then there exist constants $C_1$, $C_2>0$ such that
$$C_1\left(\frac{|E|}{|B|}\right)^p\le\frac{w(E)}{w(B)}\le C_2\left(\frac{|E|}{|B|}\right)^{(r-1)/r}$$for any measurable subset $E$ of a ball $B$.
\end{lemma}
Given a Muckenhoupt's weight function $w$ on $\mathbb R^n$, for $0<p<\infty$, we denote by $L^p(w)$ the space of all functions satisfying
\begin{equation*}
\|f\|_{L^p(w)}=\left(\int_{\mathbb R^n}|f(x)|^pw(x)\,dx\right)^{1/p}<\infty.
\end{equation*}

Throughout this article, we will use $C$ to denote a positive constant, which is independent of the main parameters and not necessarily the same at each occurrence. By $A\sim B$, we mean that there exists a constant $C>1$ such that $\frac1C\le\frac AB\le C$. Moreover, we denote the conjugate exponent of $q>1$ by $q'=q/(q-1).$
\section{Weighted Hardy spaces $H^p_L(w)$ for $0<p<1$ and their atomic decompositions}
Let $L=-\Delta+V$. For any $t>0$, we define $P_t=e^{-tL}$ and
\begin{equation*}
Q_{t,k}=(-t)^k\left.\frac{d^kP_s}{ds^k}\right|_{s=t}=(tL)^k e^{-tL},\quad k=1,2,\ldots.
\end{equation*}
We denote simply by $Q_t$ when $k=1$. First note that Gaussian upper bounds carry over from heat kernels to their time derivatives.
\begin{lemma}[{[2,14]}]
For every $k=1,2,\ldots,$ there exist two positive constants $C_k$ and $c_k$ such that the kernel $p_{t,k}(x,y)$ of the operator $Q_{t,k}$ satisfies
\begin{equation*}
\big|p_{t,k}(x,y)\big|\le\frac{C_k}{(4\pi t)^{n/2}}e^{-\frac{|x-y|^2}{c_k t}}
\end{equation*}
for all $t>0$ and almost all $x,y\in\mathbb R^n$.
\end{lemma}
Set
\begin{equation*}
H^2(\mathbb R^n)=\overline{\mathcal R(L)}=\overline{\{Lu\in L^2(\mathbb R^n):u\in L^2(\mathbb R^n)\}},
\end{equation*}
where $\overline{\mathcal R(L)}$ stands for the range of $L$. We also set
\begin{equation*}
\Gamma(x)=\{(y,t)\in{\mathbb R}^{n+1}_+:|x-y|<t\}.
\end{equation*}For a given function $f\in L^2(\mathbb R^n)$, we consider the area integral function associated to Schr\"odinger operator $L$
\begin{equation*}
S_L(f)(x)=\bigg(\iint_{\Gamma(x)}\big|Q_{t^2}(f)(y)\big|^2\frac{dydt}{t^{n+1}}\bigg)^{1/2},\quad x\in\mathbb R^n.
\end{equation*}
Given a weight function $w$ on $\mathbb R^n$, for $0<p<1$, we shall define the weighted Hardy spaces $H^p_L(w)$ as the completion of $H^2(\mathbb R^n)$ in the norm given by the $L^p(w)$-norm of area integral function; that is
\begin{equation*}
\|f\|_{H^p_L(w)}=\|S_L(f)\|_{L^p(w)}.
\end{equation*}
Let $M\in\mathbb N$ and $0<p<1$. As in [18], we say that a function $a(x)\in L^2(\mathbb R^n)$ is called a $(p,M)$-atom with respect to $w$(or a $w$-$(p,M)$-atom) if there exist a ball $B=B(x_0,r)$ and a function $b\in\mathcal D(L^M)$ such that

(a) $a=L^M b$;

(b) $supp\,L^k b\subseteq B, \quad k=0,1,\ldots,M$;

(c) $\|(r^2L)^kb\|_{L^2(B)}\le r^{2M}|B|^{1/2}w(B)^{-1/p},\quad k=0,1,\ldots,M$.

Let $M\in\mathbb N$ and $\frac{n}{n+1}<p<1$. For any $w$-$(p,M)$-atom $a$ associated to a ball $B=B(x_0,r)$, $\|a\|_{L^2(B)}\le |B|^{1/2}w(B)^{-1/p}$, we will show that $a\in H^p_L(w)$ and its $H^p_L(w)$-norm is uniformly bounded; precisely
\begin{theorem}
Let $M\in\mathbb N$, $\frac{n}{n+1}<p<1$ and $w\in A_1\cap RH_{(2/p)'}$. Then there exists a constant $C>0$ independent of $a$ such that
\begin{equation*}
\|S_L(a)\|_{L^p(w)}\le C.
\end{equation*}
\end{theorem}
\begin{proof}
We write
\begin{equation*}
\begin{split}
\big\|S_L(a)\big\|^p_{L^p(w)}&=\int_{2B}\big|S_L(a)(x)\big|^pw(x)\,dx+\int_{(2B)^c}\big|S_L(a)(x)\big|^pw(x)\,dx\\
&=I_1+I_2.
\end{split}
\end{equation*}
Set $q=2/p$. Note that $w\in RH_{q'}$, then it follows from H\"older's inequality, Lemma 2.1 and the $L^2$ boundedness of $S_L$(see (3.2) below) that
\begin{equation*}
\begin{split}
I_1&\le\Big(\int_{2B}\big|S_L(a)(x)\big|^2\,dx\Big)^{p/2}\Big(\int_{2B}w(x)^{q'}\,dx\Big)^{1/{q'}}\\
&\le C\|a\|^p_{L^2(B)}\cdot\frac{w(2B)}{|2B|^{1/q}}\\
&\le C.
\end{split}
\end{equation*}

We turn to deal with $I_2$. By using H\"older's inequality and the fact that $w\in RH_{q'}$, we can get
\begin{align*}
I_2&=\sum_{k=1}^\infty\int_{2^{k+1}B\backslash 2^k B}\big|S_L(a)(x)\big|^pw(x)\,dx\\
&\le C\sum_{k=1}^\infty\Big(\int_{2^{k+1}B\backslash 2^k B}\big|S_L(a)(x)\big|^2\,dx\Big)^{p/2}\cdot\frac{w(2^{k+1}B)}{|2^{k+1}B|^{1/q}}.
\end{align*}
For any $x\in 2^{k+1}B\backslash2^k B$, $k=1,2,\ldots,$ we write
\begin{equation*}
\begin{split}
&\big|S_L(a)(x)\big|^2\\
=&\int_0^r\int_{|y-x|<t}\big|t^2Le^{-t^2L}a(y)\big|^2\frac{dydt}{t^{n+1}}+\int_r^\infty\int_{|y-x|<t}\big|t^2Le^{-t^2L}a(y)\big|^2\frac{dydt}{t^{n+1}}\\
=&\,\mbox{\upshape I+II}.
\end{split}
\end{equation*}
For the term I, note that $0<t<r$. By a simple calculation, we obtain that for any $(y,t)\in\Gamma(x)$, $x\in 2^{k+1}B\backslash2^k B$, $z\in B$, then $|y-z|\ge2^{k-1}r$. Hence, by using H\"older's inequality and Lemma 3.1, we deduce
\begin{equation*}
\begin{split}
\big|t^2Le^{-t^2L}a(y)\big|&\le C\cdot\frac{t}{(2^{k-1} r)^{n+1}}\int_B|a(z)|\,dz\\
&\le C\cdot\frac{t}{(2^k r)^{n+1}}\|a\|_{L^2(\mathbb R^n)}|B|^{1/2}\\
&\le C\cdot w(B)^{-1/p}\frac{t}{2^{k(n+1)}\cdot r}.
\end{split}
\end{equation*}
Consequently
\begin{equation*}
\begin{split}
\mbox{\upshape I}&\le C\Big(\frac{1}{2^{k(n+1)}w(B)^{1/p}}\Big)^2\cdot\frac{1}{r^2}\int_0^r t\,dt\\
&\le C\Big(\frac{1}{2^{k(n+1)}w(B)^{1/p}}\Big)^2.
\end{split}
\end{equation*}
We now estimate the other term II. In this case, a direct computation shows that for any $(y,t)\in\Gamma(x)$, $x\in 2^{k+1}B\backslash2^k B$ and $z\in B$, we have $t+|y-z|\ge 2^{k-1}r$. Since there exists a function $b\in\mathcal D(L^M)$ such that $a=L^M b$, then by H\"older's inequality and Lemma 3.1 again, we get
\begin{equation*}
\begin{split}
\big|t^2Le^{-t^2L}a(y)\big|&=\big|(t^2L)^{M+1}e^{-t^2L}b(y)\big|\cdot\frac{1}{t^{2M}}\\
&\le C\cdot\frac{1}{(2^{k-1}r)^{n+1}}\int_B|b(z)|\,dz\cdot\frac{1}{t^{2M-1}}\\
&\le C\cdot\frac{1}{(2^kr)^{n+1}}\|b\|_{L^2(\mathbb R^n)}|B|^{1/2}\cdot\frac{1}{t^{2M-1}}\\
&\le C\cdot\frac{r^{2M-1}}{2^{k(n+1)}w(B)^{1/p}}\cdot\frac{1}{t^{2M-1}}.
\end{split}
\end{equation*}
Therefore
\begin{equation*}
\begin{split}
\mbox{\upshape II}&\le C\Big(\frac{1}{2^{k(n+1)}w(B)^{1/p}}\Big)^2\cdot r^{4M-2}\int_r^\infty\frac{dt}{t^{4M-1}}\\
&\le C\Big(\frac{1}{2^{k(n+1)}w(B)^{1/p}}\Big)^2.
\end{split}
\end{equation*}
Combining the above estimates for I and II, we thus obtain
\begin{equation*}
\big|S_L(a)(x)\big|\le C\cdot\frac{1}{2^{k(n+1)}w(B)^{1/p}},\quad \mbox{when }x\in 2^{k+1}B\backslash2^k B.
\end{equation*}
Then it follows immediately from Lemma 2.1 that
\begin{equation*}
\begin{split}
I_2&\le C\sum_{k=1}^\infty\frac{1}{2^{kp(n+1)}w(B)}\cdot w(2^{k+1}B)\\
&\le C\sum_{k=1}^\infty\frac{1}{2^{kp(n+1)-kn}}\\
&\le C,
\end{split}
\end{equation*}
where in the last inequality we have used the fact that $p>n/{(n+1)}$. Summarizing the estimates for $I_1$ and $I_2$ derived above, we complete the proof of Theorem 3.2.
\end{proof}
For every bounded Borel function $F:[0,\infty)\to\mathbb C$, we define the operator $F(L):L^2(\mathbb R^n)\to L^2(\mathbb R^n)$ by the following formula
\begin{equation*}
F(L)=\int_0^\infty F(\lambda)\,dE_L(\lambda),
\end{equation*}
where $E_L(\lambda)$ is the spectral decomposition of $L$. Therefore, the operator $\cos(t\sqrt L)$ is well-defined on $L^2(\mathbb R^n)$. Moreover, it follows from [16] that there exists a constant $c_0$ such that the Schwartz kernel $K_{\cos(t\sqrt L)}(x,y)$ of $\cos(t\sqrt L)$ has support contained in $\{(x,y)\in \mathbb R^n\times\mathbb R^n:|x-y|\le c_0t\}$. By the functional calculus for $L$ and Fourier inversion formula, whenever $F$ is an even bounded Borel function with $\hat F\in L^1(\mathbb R)$, we can write
\begin{equation*}
F(\sqrt L)=(2\pi)^{-1}\int_{-\infty}^\infty \hat F(t)\cos(t\sqrt L)\,dt.
\end{equation*}
\begin{lemma}[{[10]}]
Let $\varphi\in C^\infty_0(\mathbb R)$ be even and $supp\,\varphi\subseteq[-c_0^{-1},c_0^{-1}]$. Let $\Phi$ denote the Fourier transform of $\varphi$. Then for each $j=0,1,\ldots$, and for all $t>0$, the Schwartz kernel $K_{(t^2L)^j\Phi(t\sqrt L)}(x,y)$ of $(t^2L)^j\Phi(t\sqrt L)$ satisfies
\begin{equation*}
supp\,K_{(t^2L)^j\Phi(t\sqrt L)}\subseteq\{(x,y)\in\mathbb R^n\times\mathbb R^n:|x-y|\le t\}.
\end{equation*}
\end{lemma}
For a given $s>0$, we set
\begin{equation*}
\mathcal F(s)=\Big\{\psi:\mathbb C\to\mathbb C \mbox{ measurable}, |\psi(z)|\le C\frac{|z|^s}{1+|z|^{2s}}\Big\}.
\end{equation*}
Then for any nonzero function $\psi\in\mathcal F(s)$, we have the following estimate(see [18])
\begin{equation}
\Big(\int_0^\infty\|\psi(t\sqrt L)f\|^2_{L^2(\mathbb R^n)}\frac{dt}{t}\Big)^{1/2}\le C\|f\|_{L^2(\mathbb R^n)}.
\end{equation}
In particular, we have
\begin{equation}
\|S_L(f)\|_{L^2(\mathbb R^n)}\le C\|f\|_{L^2(\mathbb R^n)}.
\end{equation}

We are going to establish the atomic decomposition for the weighted Hardy spaces $H^p_L(w)$($0<p<1$).
\begin{theorem}
Let $M\in\mathbb N$, $0<p<1$ and $w\in A_1$. If $f\in H^p_L(w)$, then there exist a family of $w$-$(p,M)$-atoms \{$a_j$\} and a sequence of real numbers \{$\lambda_j$\} with $\sum_{j}|\lambda_j|^p\le C\|f\|^p_{H^p_L(w)}$ such that $f$ can be represented in the form $f(x)=\sum_{j}\lambda_ja_j(x)$, and the sum converges both in the sense of $L^2(\mathbb R^n)$-norm and $H^p_L(w)$-norm.
\end{theorem}
\begin{proof}
First assume that $f\in H^p_L(w)\cap H^2(\mathbb R^n)$. We follow the same constructions as in [18]. Let $\varphi$ and $\Phi$ be as in Lemma 3.3. We set $\Psi(x)=x^{2M}\Phi(x), x\in\mathbb R$. By the $L^2$-functional calculus of $L$, for every $f\in H^2(\mathbb R^n)$, we can establish the following version of the Calder\'on reproducing formula
\begin{equation}
f(x)=c_\psi\int_0^\infty\Psi(t\sqrt L)t^2 Le^{-t^2 L}(f)(x)\frac{dt}{t},
\end{equation}
where the above equality holds in the sense of $L^2(\mathbb R^n)$-norm. For any $k\in\mathbb Z$, set
$$\Omega_k=\{x\in\mathbb R^n:S_{L,10\sqrt n}(f)(x)>2^k\}.$$
Let $\mathbb D$ denote the set formed by all dyadic cubes in $\mathbb R^n$ and let
\begin{equation*}
\mathbb D_k=\Big\{Q\in{\mathbb D}:w(Q\cap\Omega_k)>\frac{w(Q)}{2},w(Q\cap\Omega_{k+1})\le\frac{w(Q)}{2}\Big\}.
\end{equation*}
Obviously, for any $Q\in\mathbb D$, there exists a unique $k\in\mathbb Z$ such that $Q\in{\mathbb D}_k.$
We also denote the maximal dyadic cubes in ${\mathbb D}_k$ by $Q_k^l$. Note that the maximal dyadic cubes $Q_k^l$ are pairwise disjoint, then it is easy to check that
\begin{equation}
\sum_lw(Q^l_k)\le C\cdot w(\Omega_k).
\end{equation}Set
\begin{equation*}
\widetilde Q=\Big\{(y,t)\in{\mathbb R}^{n+1}_+:y\in Q,\, \frac{l(Q)}{2}<t\le l(Q)\Big\},
\end{equation*}
where $l(Q)$ denotes the side length of $Q$. If we set $\widetilde{Q_k^l}=\underset{Q_k^l\supseteq Q\in{\mathbb D}_k}{\bigcup}\widetilde Q$,
then we have ${\mathbb R}^{n+1}_+=\underset{k}{\bigcup}\,\underset{l}{\bigcup}\,\widetilde{Q^l_k}$. Hence, by the formula (3.3), we can write
\begin{equation*}
f(x)=\sum_k\sum_l c_\psi\int_{\widetilde{Q^l_k}}\Psi(t\sqrt L)(x,y)t^2Le^{-t^2L}f(y)\frac{dydt}{t}=\sum_k\sum_l\lambda_{kl}a^l_k(x),
\end{equation*}
where $a^l_k=L^M b^l_k$,
\begin{equation*}b^l_k(x)=c_\psi\lambda_{kl}^{-1}\int_{\widetilde{Q^l_k}}t^{2M}\Phi(t\sqrt L)(x,y)t^2Le^{-t^2L}f(y)\frac{dydt}{t}\end{equation*}
and
\begin{equation*}
\lambda_{kl}=w(Q^l_k)^{1/p-1/2}\bigg(\int_{\widetilde{Q^l_k}}\big|t^2Le^{-t^2L}f(y)\big|^2\frac{w(Q^l_k)}{|Q^l_k|}\frac{dydt}{t}\bigg)^{1/2}.
\end{equation*}
By using Lemma 3.3, the authors in [18] showed that for every $j=0,1,\ldots,M$, $supp\,(L^j b^l_k)\subseteq3Q^l_k$. In [18], they also obtained the following estimate
\begin{equation}
\sum_l\int_{\widetilde{Q^l_k}}\big|t^2Le^{-t^2L}f(y)\big|^2\frac{w(Q^l_k)}{|Q^l_k|}\frac{dydt}{t}\le C\cdot2^{2k}w(\Omega_k).
\end{equation}
Since
\begin{equation*}
\big\|\big(l(Q^l_k)^2L\big)^jb^l_k\big\|_{L^2(3Q^l_k)}=\sup_{\|h\|_{L^2(3Q^l_k)}\le1}\bigg|\int_{\mathbb R^n}\big(l(Q^l_k)^2L\big)^jb^l_k(x)h(x)\,dx\bigg|.
\end{equation*}
Let $\Psi_j(x)=x^{2j}\Phi(x)$, $j=0,1,\ldots, M$. Then we can easily verify that $\Psi_j\in\mathcal F(2j)$. Observe that when $(y,t)\in\widetilde{Q^l_k}$, we have $t\sim l(Q^l_k)$. Then it follows from H\"older's inequality and the estimate (3.1) that
\begin{equation*}
\begin{split}
&\bigg|\int_{\mathbb R^n}\big(l(Q^l_k)^2L\big)^jb^l_k(x)h(x)\,dx\bigg|\\
\le& \frac{C\cdot l(Q^l_k)^{2M}}{\lambda_{kl}}\bigg(\int_{\widetilde{Q^l_k}}\big|t^2Le^{-t^2L}f(y)\big|^2\frac{dydt}{t}\bigg)^{1/2}\\
&\times\bigg(\int_{\widetilde{Q^l_k}}\Big|(t^2L)^j\Phi(t\sqrt L)(h\chi_{3Q^l_k})(y)\Big|^2\frac{dydt}{t}\bigg)^{1/2}\\
\le& C\cdot l(Q^l_k)^{2M}|Q^l_k|^{1/2}w(Q^l_k)^{-1/p}\Big(\int_0^\infty\|\Psi_j(t\sqrt L)(h\chi_{3Q^l_k})\|^2_{L^2(\mathbb R^n)}\frac{dt}{t}\Big)^{1/2}\\
\le& C\cdot l(Q^l_k)^{2M}|Q^l_k|^{1/2}w(Q^l_k)^{-1/p}.
\end{split}
\end{equation*}
Hence
\begin{equation*}
\big\|\big(l(Q^l_k)^2L\big)^jb^l_k\big\|_{L^2(3Q^l_k)}\le C\cdot l(Q^l_k)^{2M}|Q^l_k|^{1/2}w(Q^l_k)^{-1/p}.
\end{equation*}
From the above discussions, we have proved that these functions $a^l_k$ are all $w$-$(p,M)$-atoms up to a normalization by a multiplicative constant. Finally, by using H\"older's inequality, the estimates (3.4) and (3.5), we obtain
\begin{equation*}
\begin{split}
\sum_k\sum_l|\lambda_{kl}|^p&=\sum_k\sum_l\Big(w(Q^l_k)\Big)^{1-p/2}\bigg(\int_{\widetilde{Q^l_k}}\big|t^2Le^{-t^2L}f(y)\big|^2\frac{w(Q^l_k)}{|Q^l_k|}\frac{dydt}{t}\bigg)^{p/2}\\
&\le\sum_k\Big(\sum_lw(Q^l_k)\Big)^{1-p/2}\bigg(\sum_l\int_{\widetilde{Q^l_k}}\big|t^2Le^{-t^2L}f(y)\big|^2\frac{w(Q^l_k)}{|Q^l_k|}\frac{dydt}{t}\bigg)^{p/2}\\
&\le C\sum_k\Big(w(\Omega_k)\Big)^{1-p/2}\Big(2^{2k}w(\Omega_k)\Big)^{p/2}\\
&\le C\|S_L(f)\|^p_{L^p(w)}.
\end{split}
\end{equation*}
Therefore, we have established the atomic decomposition for all functions in the space $H^p_L(w)\cap H^2(\mathbb R^n)$. By a standard density argument, we can show that the same conclusion holds for $H^p_L(w)$. Following along the same arguments as in [18], we can also prove that the sum $f=\sum_{j}\lambda_ja_j$ converges both in the sense of $L^2(\mathbb R^n)$-norm and $H^p_L(w)$-norm, the details are omitted here. This completes the proof of Theorem 3.4.
\end{proof}

\section{Proof of Theorem 1.1}
We shall need the following Davies-Gaffney estimate which can be found in [10,18].
\begin{lemma}
For any two closed sets $E$ and $F$ of $\mathbb R^n$, there exist two positive constants $C$ and $c$ such that
\begin{equation*}
\|t\nabla e^{-t^2L}f\|_{L^2(F)}\le C \cdot e^{-\frac{dist(E,F)^2}{ct^2}}\|f\|_{L^2(E)}
\end{equation*}
for every $f\in L^2(\mathbb R^n)$ with support contained in $E$.
\end{lemma}
\begin{theorem}
Let $\frac{n}{n+1}<p<1$ and $w\in A_1\cap RH_{(2/p)'}$. Then the operator $\nabla L^{-1/2}$ is bounded from $H^p_L(w)$ to $L^p(w)$.
\end{theorem}
\begin{proof}
By Theorem 3.4 we just proved, it is enough for us to show that for any $w$-$(p,M)$-atom $a$, $M>\frac n2(\frac1p-\frac12)$, there exists a constant $C>0$ independent of $a$ such that
$\|\nabla L^{-1/2}(a)\|_{L^p(w)}\le C$. Let $a$ be a $w$-$(p,M)$-atom with $supp\,a\subseteq B=B(x_0,r)$, $\|a\|_{L^2(B)}\le |B|^{1/2}w(B)^{-1/p}$. We write
\begin{equation*}
\begin{split}
\big\|\nabla L^{-1/2}(a)\big\|^p_{L^p(w)}&=\int_{2B}\big|\nabla L^{-1/2}(a)(x)\big|^pw(x)\,dx+\int_{(2B)^c}\big|\nabla L^{-1/2}(a)(x)\big|^pw(x)\,dx\\
&=J_1+J_2.
\end{split}
\end{equation*}
Set $q=2/p$. Applying H\"older's inequality, the $L^2$ boundedness of $\nabla L^{-1/2}$, Lemma 2.1 and $w\in RH_{q'}$, we thus have
\begin{equation*}
\begin{split}
J_1&\le\Big(\int_{2B}\big|\nabla L^{-1/2}(a)(x)\big|^2\,dx\Big)^{p/2}\Big(\int_{2B}w(x)^{q'}\,dx\Big)^{1/{q'}}\\
&\le C\|a\|^p_{L^2(B)}\cdot\frac{w(2B)}{|2B|^{1/q}}\\
&\le C.
\end{split}
\end{equation*}
On the other hand, it follows from H\"older's inequality and $w\in RH_{q'}$ that
\begin{align}
J_2&=\sum_{k=1}^\infty\int_{2^{k+1}B\backslash 2^k B}\big|\nabla L^{-1/2}(a)(x)\big|^pw(x)\,dx\notag\\
&\le C\sum_{k=1}^\infty\Big(\int_{2^{k+1}B\backslash 2^k B}\big|\nabla L^{-1/2}(a)(x)\big|^2\,dx\Big)^{p/2}\cdot\frac{w(2^{k+1}B)}{|2^{k+1}B|^{1/q}}.
\end{align}
By a change of variable $s=t^2$, we can rewrite (1.2) as
\begin{equation}
\nabla L^{-1/2}(a)(x)=\frac{2}{\sqrt\pi}\int_0^\infty s\nabla e^{-{s}^2L}(a)(x)\frac{ds}{s}.
\end{equation}
For any $k=1,2,\ldots,$ it follows immediately from Minkowski's integral inequality that
\begin{equation*}
\begin{split}
&\Big(\int_{2^{k+1}B\backslash 2^k B}\big|\nabla L^{-1/2}(a)(x)\big|^2\,dx\Big)^{1/2}\\
\le& \,C\int_0^r\big\|s\nabla e^{-{s}^2L}a\big\|_{L^2(2^{k+1}B\backslash 2^k B)}\frac{ds}{s}+C\int_r^\infty\big\|s\nabla e^{-{s}^2L}a\big\|_{L^2(2^{k+1}B\backslash 2^k B)}\frac{ds}{s}\\
=&\,\mbox{\upshape III+IV}.
\end{split}
\end{equation*}
Observe that $M>\frac n2(\frac1p-\frac12)$. Then we are able to choose a positive number $N$ such that $\frac n2(\frac1p-\frac12)<N<M$. By using Lemma 4.1, we can get
\begin{align}
\mbox{\upshape III}&\le C\int_0^r e^{-\frac{(2^kr)^2}{s^2}}\|a\|_{L^2(B)}\frac{ds}{s}\notag\\
&\le C\int_0^r\frac{s^{2N}}{(2^kr)^{2N}}\frac{ds}{s}\cdot\|a\|_{L^2(B)}\notag\\
&\le C\cdot2^{-2kN}|B|^{1/2}w(B)^{-1/p}.
\end{align}
We now turn to estimate the term IV. Since $a=L^M b$ and $\|b\|_{L^2(B)}\le r^{2M}|B|^{1/2}w(B)^{-1/p}$. Using Lemma 3.1 and Lemma 4.1, we deduce
\begin{align}
\mbox{\upshape IV}&=C\int_r^\infty\big\|s\nabla e^{-{s}^2L}(L^Mb)\big\|_{L^2(2^{k+1}B\backslash 2^k B)}\frac{ds}{s}\notag\\
&=C\int_r^\infty\big\|s\nabla e^{-\frac{s^2L}{2}}(s^2L)^Me^{-\frac{s^2L}{2}}b\big\|_{L^2(2^{k+1}B\backslash 2^k B)}\frac{ds}{s^{2M+1}}\notag\\
&\le C\int_r^\infty e^{-\frac{(2^kr)^2}{s^2}}\big\|(s^2L)^Me^{-\frac{s^2L}{2}}b\big\|_{L^2(B)}\frac{ds}{s^{2M+1}}\notag\\
&\le C\int_r^\infty\frac{s^{2N}}{(2^kr)^{2N}}\frac{ds}{s^{2M+1}}\cdot\|b\|_{L^2(B)}\notag\\
&\le C\cdot2^{-2kN}|B|^{1/2}w(B)^{-1/p}.
\end{align}
Combining the above inequality (4.4) with (4.3), we thus obtain
\begin{equation}
\Big(\int_{2^{k+1}B\backslash 2^k B}\big|\nabla L^{-1/2}(a)(x)\big|^2\,dx\Big)^{1/2}\le C\cdot2^{-2kN}|B|^{1/2}w(B)^{-1/p}.
\end{equation}
Substituting the above inequality (4.5) into (4,1) and using Lemma 2.1, then we have
\begin{equation*}
\begin{split}
J_2&\le C\sum_{k=1}^\infty\Big(2^{-2kN}|B|^{1/2}w(B)^{-1/p}\Big)^p\cdot\frac{w(2^{k+1}B)}{|2^{k+1}B|^{p/2}}\\
&\le C\sum_{k=1}^\infty\frac{1}{2^{k(2pN+\frac{np}{2}-n)}}\\
&\le C,
\end{split}
\end{equation*}
where the last series is convergent since $N>\frac{n}{2}(\frac1p-\frac12)$. Summarizing the estimates for $J_1$ and $J_2$, we get the desired result.
\end{proof}
The real-variable theory of classical weighted Hardy spaces have been extensively studied by many authors. In 1979, Garcia-Cuerva studied the atomic decomposition and the dual spaces of $H^p(w)$ for $0<p\le1$. In 2002, Lee and Lin gave the molecular characterization of $H^p(w)$ for $0<p\le1$, they also obtained the $H^p(w)$($\frac12<p\le1$) boundedness of the Hilbert transform and the $H^p(w)$($\frac n{n+1}<p\le1$) boundedness of the Riesz transforms. For the results mentioned above, we refer the readers to [7,12,19] for further details.

Let $\frac{n}{n+1}<p<1$ and $w\in A_1$. A real-valued function $a(x)$ is called a $w$-$(p,2,0)$-atom if the following conditions are satisfied(see [7,19]):

(a) $supp\,a\subseteq B$;

(b) $\|a\|_{L^2(B)}\le |B|^{1/2}w(B)^{-1/p}$;

(c) $\int_{\mathbb R^n}a(x)\,dx=0$.
\begin{theorem}
Let $\frac{n}{n+1}<p<1$ and $w\in A_1$. For each $f\in H^p(w)$, there exist a family of $w$-$(p,2,0)$-atoms \{$a_j$\} and a sequence of real numbers \{$\lambda_j$\} with $\sum_j|\lambda_j|^p\le C\|f\|^p_{H^p(w)}$ such that $f=\sum_j\lambda_j a_j$ in the sense of $H^p(w)$ norm.
\end{theorem}
Next, as in [18], we shall also define the new weighted molecules for $H^p(w)$. Let $\frac{n}{n+1}<p<1$, $w\in A_1$ and $\varepsilon>0$. A function $m(x)\in L^2(\mathbb R^n)$ is called a $w$-$(p,2,0,\varepsilon)$-molecule associated to a ball $B$ if the following conditions are satisfied:

(A) $\int_{\mathbb R^n}m(x)\,dx=0$;

(B) $\|m\|_{L^2(2B)}\le |B|^{1/2}w(B)^{-1/p}$;

(C) $\|m\|_{L^2(2^{k+1}B\backslash 2^k B)}\le 2^{-k\varepsilon}|2^kB|^{1/2}w(2^kB)^{-1/p},\quad k=1,2,\ldots$.

Note that for every $w$-$(p,2,0)$-atom, it is a $w$-$(p,2,0,\varepsilon)$-molecule for all $\varepsilon>0$. Then we are able to establish the following molecular characterization for the classical weighted Hardy spaces $H^p(w)$.
\begin{theorem}
Let $\frac{n}{n+1}<p<1$ and $w\in A_1$.

\noindent$(i)$ If $f\in H^p(w)$, then there exist a family of $w$-$(p,2,0,\varepsilon)$-molecules \{$m_j$\} and a sequence of real numbers \{$\lambda_j$\} with $\sum_j|\lambda_j|^p\le C\|f\|^p_{H^p(w)}$ such that $f=\sum_j\lambda_j a_j$ in the sense of $H^p(w)$ norm.

\noindent$(ii)$ Suppose that $w\in RH_{(2/p)'}$ and $\varepsilon>n/2$, then every $w$-$(p,2,0,\varepsilon)$-molecule $m$ is in $H^p(w)$. Moreover, there exists a constant $C>0$ independent of $m$ such that $\|m\|_{H^p(w)}\le C$.
\end{theorem}
\begin{proof}
$(i)$ is a straightforward consequence of Theorem 4.3.

$(ii)$ We follow the idea of [18]. Denote $m_0(x)=m(x)\chi_{_{2B}}(x)$, $m_k(x)=m(x)\chi_{_{2^{k+1}B\backslash 2^k B}}(x)$, $k=1,2,\ldots.$ Then we can decompose $m(x)$ as
\begin{equation*}
m(x)=\sum_{k=0}^\infty m_k(x)=\sum_{k=0}^\infty\big(m_k(x)-N_k(x)\big)+\sum_{k=0}^\infty N_k(x),
\end{equation*}
where $N_0(x)=\frac{1}{|2B|}\int_{\mathbb R^n}m_0(y)\,dy\cdot\chi_{_{2B}}(x)$ and $N_k(x)=\frac{1}{|2^{k+1}B\backslash 2^k B|}\int_{\mathbb R^n}m_k(y)\,dy\cdot\chi_{_{2^{k+1}B\backslash 2^k B}}(x)$, $k=1,2,\ldots.$ Following along the same lines as in [18], we can also show that each $(m_k-N_k)$ is a multiple of $w$-$(p,2,0)$-atom with a sequence of coefficients in $\mathit l^p$. We set $\eta_k=\int_{\mathbb R^n}m_k(y)\,dy$, $k=0,1,\ldots.$ In [18], Song and Yan established the following identity
\begin{equation*}
\sum_{k=0}^\infty N_k(x)=\sum_{k=0}^\infty p_k\cdot\psi_k(x),
\end{equation*}
where $p_k=\sum_{j=k+1}^\infty\eta_j$ and $\psi_k(x)=\frac{N_{k+1}(x)}{\eta_{k+1}}-\frac{N_k(x)}{\eta_k}$. Then we have
\begin{equation*}
\begin{split}
\big|p_k\big|&\le \sum_{j=k+1}^\infty\int_{2^{j+1}B\backslash 2^j B}|m(y)|\,dy\\
&\le\sum_{j=k+1}^\infty\|m\|_{L^2(2^{j+1}B\backslash 2^j B)}|2^{j+1}B|^{1/2}\\
&\le C\sum_{j=k+1}^\infty2^{-j\varepsilon}\cdot|2^jB|w(2^j B)^{-1/p}.
\end{split}
\end{equation*}
When $j\ge k+1$, then $2^k B\subseteq2^j B$. Since $w\in RH_{(2/p)'}$, then by Lemma 2.2, we can get
\begin{equation*}
\frac{w(2^k B)}{w(2^j B)}\le C\left(\frac{|2^k B|}{|2^j B|}\right)^{p/2}.
\end{equation*}
Hence
\begin{equation*}
\begin{split}
\big|p_k\big|&\le C\cdot\frac{|2^k B|}{w(2^kB)^{1/p}}\sum_{j=k+1}^\infty2^{-j\varepsilon}\left(\frac{|2^jB|}{|2^kB|}\right)^{1/2}\\
&\le C\cdot\frac{|2^k B|}{w(2^kB)^{1/p}}\bigg(\sum_{j=k+1}^\infty2^{-j(\varepsilon-n/2)}\bigg)\cdot2^{-kn/2}\\
&\le C\cdot2^{-k\varepsilon}\frac{|2^k B|}{w(2^kB)^{1/p}},
\end{split}
\end{equation*}
where the last inequality holds since $\varepsilon>n/2$. As in [18], we can easily check that $2^{k\varepsilon}p_k\psi(x)$ are all $w$-$(p,2,0)$-atoms associated to $2^{k+1}B$. Therefore the sum $\sum_{k=0}^\infty N_k$ can be write as an infinite linear combination of $w$-$(p,2,0)$-atoms with a sequence of coefficients in $\mathit l^p$. Summarizing the above discussions, we complete the proof of Theorem 4.4.
\end{proof}
We are now in a position to give the proof of Theorem 1.1.
\begin{proof}[Proof of Theorem 1.1]
By Theorem 3.4 and Theorem 4.4, it suffices to show that for every $w$-$(p,M)$-atom $a$ with $supp\, a\subseteq B$, then $\nabla L^{-1/2}a$ is a $w$-$(p,2,0,\varepsilon)$-molecule, where $M>n(\frac1p-\frac12)$ and $\varepsilon>n/2$. It is easy to see that $\int_{\mathbb R^n}\nabla L^{-1/2}a(x)\,dx=0$. It remains to verify the estimates (B) and (C). H\"older's inequality and the definition of $w$-$(p,M)$-atom imply
\begin{equation*}
\big\|\nabla L^{-1/2}(a)\big\|_{L^2(2B)}\le C\|a\|_{L^2(B)}\le C\cdot|B|^{1/2}w(B)^{-1/p}.
\end{equation*}
For $k=1,2,\ldots$, it follows from the previous estimates (4.3) and (4.4) that
\begin{equation*}
\big\|\nabla L^{-1/2}(a)\big\|_{L^2(2^{k+1}B\backslash 2^k B)}\le C\cdot2^{-2kN}|B|^{1/2}w(B)^{-1/p},
\end{equation*}
where $N>0$ is chosen such that $n(\frac1p-\frac12)<N<M$. By using Lemma 2.2, we get
\begin{equation*}
\frac{w(B)}{w(2^kB)}\ge C\cdot\frac{|B|}{|2^k B|}.
\end{equation*}
Hence
\begin{equation*}
\big\|\nabla L^{-1/2}(a)\big\|_{L^2(2^{k+1}B\backslash 2^k B)}\le C\cdot2^{-k(2N-n/p+n/2)}|2^kB|^{1/2}w(2^kB)^{-1/p}.
\end{equation*}
Therefore, we have proved $\nabla L^{-1/2}a$ is a $w$-$(p,2,0,2N-n/p+n/2)$-molecule. This concludes the proof of Theorem 1.1. 
\end{proof}

\end{document}